\documentclass[a4paper, reqno, 11pt]{amsart}

\usepackage{amscd, amssymb, amsthm, amsmath, color}
\usepackage{graphicx, psfrag, enumerate}
\newtheorem{lem}{Lemma}[section]
\newtheorem{thm}[lem]{Theorem}
\newtheorem{pro}[lem]{Proposition}
\newtheorem{cor}[lem]{Corollary}

\newtheorem{defi}[lem]{Definition}

\newcommand{\Sor}{{\mathsf{Sor}}}
\newcommand{\sor}{{\mathsf{sor}}}
\newcommand{\maj}{{\mathsf{maj}}}
\newcommand{\Des}{{\mathsf{Des}}}

\newcommand{\Inv}{{\mathsf{Inv}}}
\newcommand{\inv}{{\mathsf{inv}}}

\newcommand{\Leh}{{\mathsf{Leh}}}

\newcommand{\stat}{{\mathsf{stat}}}
\newcommand{\Stat}{{\mathsf{Stat}}}
\newcommand{\Rmap}{{\mathsf{Rmap}}}
\newcommand{\Rmal}{{\mathsf{Rmal}}}
\newcommand{\rlmax}{{\mathsf{rmax}}}
\newcommand{\Rmip}{{\mathsf{Rmip}}}
\newcommand{\Rmil}{{\mathsf{Rmil}}}
\newcommand{\rlmin}{{\mathsf{rmin}}}
\newcommand{\Lmap}{{\mathsf{Lmap}}}
\newcommand{\Lmal}{{\mathsf{Lmal}}}
\newcommand{\lrmax}{{\mathsf{lmax}}}
\newcommand{\Lmip}{{\mathsf{Lmip}}}
\newcommand{\Lmil}{{\mathsf{Lmil}}}
\newcommand{\lrmin}{{\mathsf{lmin}}}

\newcommand{\Lmicycl}{{\mathsf{Lmic}}}
\newcommand{\lrmincyc}{{\mathsf{lmic}}}

\newcommand{\Cyc}{{\mathsf{Cyc}}}
\newcommand{\cyc}{{\mathsf{cyc}}}

\newcommand{\imaj}{{\mathsf{imaj}}}
\newcommand{\rmaj}{{\mathsf{rmaj}}}
\newcommand{\chg}{{\mathsf{chg}}}
\newcommand{\cochg}{{\mathsf{cochg}}}

\begin{document}
%%%%%%%%%%%%%%%%%%%%%%%%%%%%%%%%%%%%%%%%%%%%%%%%%%%%%%%%%%%
\title{Set-valued sorting index and joint equidistributions}

\author{Sen-Peng Eu}
\address{Department of Mathematics \\
National Taiwan Normal University \\
Taipei 11677, Taiwan, ROC}
\email[Sen-Peng Eu]{speu@ntnu.edu.tw}

\author{Yuan-Hsun Lo}
%\address{Department of Mathematics \\
%National Taiwan Normal University \\
%Taipei 11677, Taiwan, ROC}
\email[Yuan-Hsun Lo]{yhlo0830@gmail.com}

\author{Tsai-Lien Wong}
\address{Department of Applied Mathematics \\
National Sun Yat-sen University \\
Kaohsiung 80424, Taiwan, ROC}
\email[Tsai-Lien Wong]{tlwong@math.nsysu.edu.tw}



\thanks{Partially supported by National Science Council, Taiwan under grants NSC 101-2115-M-003-013-MY3 (S.-P. Eu), 102-2811-M-003-025 (Y.-H. Lo) and NSC 102-2115-M-110-006-MY2 (T.-L Wong).}

%\date{\today}

\maketitle
%%%%%%%%%%%%%%%%%%%%%%%%%%%%%%%%%%%%%%%%%%%%%%%%%%%%%%%%%%%

%%%%%%%%%%%%%%%%%%%%%%%%%%%%%%%%%%%%%%%%%%%%%%%%%%%%%%%%%%%

\begin{abstract}
Recently Petersen defined a new Mahonian index $\sor$ over the symmetric group $\mathfrak{S}_n$ and proved that $(\inv, \rlmin)$ and $(\sor, \cyc)$ have the same joint distribution.
Foata and Han proved that the pairs of set-valued statistics $(\Cyc, \Rmil),(\Cyc, \Lmap), (\Rmil, \Lmap)$ have the same joint distribution over $\mathfrak{S}_n$.

In this paper we introduce the set-valued statistics $\Inv$, $\Lmil$, $\Sor$ and $\Lmicycl_1$ and generalize simultaneously
results of Petersen and Foata-Han and find many equidistributed triples of set-valued statistics and quadruples of statistics.
\end{abstract}

%%%%%%%%%%%%%%%%%%%%%%%%%%%%%%%%%%%%%%%%%%%%%%%%%%%%%%%%%%%%%%%%%%%%%%%%%%
%%%%%%%%%%%%%%%%%%%%%%%%%%%%%%%%%%%%%%%%%%%%%%%%%%%%%%%%%%%%%%%%%%%%%%%%%%

\section{Introduction} \label{sec:intro_new}
Let $\mathfrak{S}_n$ be the symmetric group of $[n]:=\{1,2,\dots ,n\}$.
For $\sigma=\sigma_1\sigma_2\dots \sigma_n \in \mathfrak{S}_n$, define the \emph{inversion} statistic $\inv$ by
$$\inv(\sigma):=\# \{(i,j): i<j \text{ and }\sigma_i<\sigma_j\},$$
and the \emph{cycle} statistic $\cyc(\sigma)$ by
$$\cyc(\sigma):= \text{the number of cycles in the cycle decomposition of $\sigma$}.$$
We say a permutation statistic over $\mathfrak{S}_n$ is \emph{Mahonian} if it is equidistributed with $\inv$, and is \emph{Stirling} if with $\cyc$.
For example, it is well known that the \emph{right to left minimum} statistic $\rlmin$, defined by
$$\rlmin(\sigma):=\# \{\sigma_i: \sigma_i < \sigma_j \text{ for all } j>i \},$$
is Stirling.

Recently Petersen found a new Mahonian statistic $\sor$, called the \emph{sorting index} (see Section 2 for definition), and proved the following:

\begin{thm}[Petersen~\cite{Petersen_11}] \label{thm:Petersen}
The pairs of statistics $(\inv, \rlmin)$ and $(\sor, \cyc)$ have the same joint distribution over $\mathfrak{S}_n$.
Also, we have
\begin{eqnarray*}
\sum_{\sigma\in\mathfrak{S}_n} q^{\inv(\sigma)}x^{\rlmin(\sigma)} &=&
\sum_{\sigma\in\mathfrak{S}_n} q^{\sor(\sigma)}x^{\cyc(\sigma)} \\ &=&
x\prod_{r=2}^{n} (x+[r]_q-1),
\end{eqnarray*}
where $[r]_q:=1+q+q^2+\dots +q^{r-1}$.
\end{thm}

A combinatorial proof is found by Chen et al.~\cite{Chen_12} via a bijection $\phi:\mathfrak{S}_n \to \mathfrak{S}_n$ such that
$$ (\inv, \rlmin) \sigma = (\sor, \cyc) \phi(\sigma),$$
where $(\inv,\rlmin) \sigma$ means $(\inv(\sigma), \rlmin(\sigma))$.
The bijection $\phi$ turns out to be of the form
$$\phi=\text{$B(\sigma)$}^{-1}\circ \text{$A(\sigma)$},$$
a composition of $B$-code and $A$-code introduced by Foata and Han~\cite{Foata_09}.
By defining the set-values statistics \emph{right-to-left minimum letters} $\Rmil$, \emph{left-to-right maximum places} $\Lmap$ and the \emph{cycle set} $\Cyc$ respectively by
$$\Rmil(\sigma):= \{ \sigma_i:
\sigma_i<\sigma_j \text{ for all } j > i\},$$
$$\Lmap(\sigma):= \{i:
\sigma_i > \sigma_j \text{ for all } i>j \},$$ and
$$\Cyc(\sigma):=\{\text{the smallest number in each cycle of the cycle
decomposition}\},$$
Foata and Han derived the following set-valued joint equidistribution results.

\begin{thm}[Foata, Han~\cite{Foata_09}] \label{thm:Foata_Han}
The followings hold.
\begin{enumerate}
\item For $\sigma\in \mathfrak{S}_n$, we have
    $$(\Rmil,\Lmap) \sigma = (\Cyc,\Lmap) \phi(\sigma).$$
\item The set-valued statistics $(\Cyc,\Rmil)$, $(\Cyc,\Lmap)$, $(\Rmil,\Lmap)$ are symmetric and joint equidistributed over $\mathfrak{S}_n$.
\end{enumerate}
\end{thm}

\medskip
The motivation of this work is to generalize above two theorems, to see if there is a set-valued version of Petersen's result or other pairs of set-valued statistics having the same distribution \`a la Foata and Han.
It turns out that we can have them both.

Throughout the paper a statistic is set-valued if and only if the first letter is in capital.
By introducing new set-valued statistics $\Inv$, $\Lmil$, $\Sor$ and $\Lmicycl_1$ (see Section 2 for definitions) and corresponding ordinary statistics $\lrmin$, $\lrmincyc_1$, our first main theorem extends simultaneously both Petersen and Foata-Han's results.

\begin{thm} \label{thm:main}
We have:
\begin{enumerate}
\item For $\sigma \in\mathfrak{S}_n$, the following holds:
$$(\Inv,\Rmil,\Lmap,\Lmil)\sigma=(\Sor,\Cyc,\Lmap,\Lmicycl_1)\phi(\sigma).$$
\item The quadruple statistics $(\inv, \rlmin, \lrmax, \lrmin)$ and
$(\sor, \cyc, \lrmax, \lrmincyc_1)$ have the same joint
distribution over $\mathfrak{S}_n$, and
\begin{eqnarray*}
\sum_{\sigma\in\mathfrak{S}_n} q^{\inv(\sigma)}x^{\rlmin(\sigma)} y^{\lrmin(\sigma)} &=&
\sum_{\sigma\in\mathfrak{S}_n} q^{\sor(\sigma)}x^{\cyc(\sigma)} y^{\lrmincyc_1(\sigma)} \\
&=& xy\prod_{r=2}^{n} (x+[r]_q+yq^{r-1}-1-q^{r-1}).
\end{eqnarray*}
\end{enumerate}
\end{thm}

Theorem~\ref{thm:main} generalizes Theorem~\ref{thm:Petersen} and Theorem~\ref{thm:Foata_Han}~(1).
We will see that Theorem~\ref{thm:Foata_Han}~(2) will be generalized in a later theorem.

By abuse of terminology, a set-valued statistic is called Mahonian
(or Stirling) if the corresponding ordinary statistic is so. Our
second result is to find triples of set-valued Stirling statistics
having the same joint distribution as $(\Rmil,\Lmap,\Lmil)$ and
$(\Cyc, \Lmicycl_1, \Lmap)$. By switching between certain set-valued
statistics by applying invese, reverse, or complement operations on
permutations, we obtain $8$ more (and $4$ partial) triples of
set-valued statistics having the same joint distribution. See
Theorem~\ref{thm:set_all} for the complete list. From these we
obtain two sets of symmetric and equidistributed pairs of set-valued
statistics, the first of which includes those three pairs in
Theorem~\ref{thm:Foata_Han}~(2).

The third part of the work is to derive quadruples of statistics
which are joint equidistributed with $(\inv, \rlmin, \lrmax,
\lrmin)$ and $(\sor, \cyc, \lrmincyc_1, \lrmax)$. Note that the
first statistic is Mahonian and the others are Stirling. Again by
switching among statistics, in Theorem~\ref{thm:sta_all} we derive
$10$ more (and $11$ partial) quadruples of statistics having the
same joint distribution.

\medskip
The rest of the paper is organized as follows.
Definitions and preliminary results will be put in Section~\ref{sec:preli}.
In Section~\ref{sec:main_proof} we prove Theorem~\ref{thm:main}.
Section~\ref{sec:triples} is devoted to triples of set-valued statistics, and Section~\ref{sec:quadruples} to quadruples of ordinary statistics.

%%%%%%%%%%%%%%%%%%%%%%%%%%%%%%%%%%%%%%%%%%%%%%%%%%%%%%%%%%%%%%%%%%%%%%%%%%
%%%%%%%%%%%%%%%%%%%%%%%%%%%%%%%%%%%%%%%%%%%%%%%%%%%%%%%%%%%%%%%%%%%%%%%%%%
\section{Preliminary results}\label{sec:preli}

\subsection{$A$-code and $B$-code}\label{sec:preli_ABcode}
We first introduce the $A$- and $B$-code of Foata and Han~\cite{Foata_09}, which are the key tools of this paper.
Given $\sigma\in \mathfrak{S}_n$, define its \emph{Lehmer code}~\cite{Lehmer_60} by
$$\Leh(\sigma):=(\ell_1,\ell_2,\ldots, \ell_n),$$
where $\ell_i=|\{j:\,1\leq j\leq i,\sigma_j\leq\sigma_i\}|.$
Let $\mathcal{L}_n:= \{(\ell_1,\ell_2,\ldots, \ell_n): 1\leq \ell_i\leq i \text{ for } 1\le i \le n\}$.
It is clear that $\Leh:\mathfrak{S}_n\to\mathcal{L}_n$ is a bijection.
The \emph{$A$-code} of a permutation $\sigma$ is defined by
$$A(\sigma):=\Leh(\sigma^{-1}).$$
For example, let $\sigma=2413765$.
Then $\sigma^{-1}=3142765$ and $A(\sigma)=(1,1,3,2,5,5,5)$.

The \emph{$B$-code} of $\sigma$ is defined in the following way.
For each $i=1,\ldots,n$, let $k_i\ge 1$ be the smallest integer such that $\sigma^{-k_i}(i)\leq i$.
Define
$$B(\sigma):=(b_1,b_2,\ldots,b_n) \text{ with } b_i=\sigma^{-k_i}(i).$$
Equivalently, $B(\sigma)$ can be determined from the cycle decomposition of $\sigma$.
Assume that $i$ appears in a cycle $c$.
If $i$ is the smallest element of $c$, then set $b_i=i$; otherwise, choose $b_i$ to be the first element $j$ in $c$
with respect to the reverse direction such that $j<i$.
For example, let $\sigma=2431756=(124)(3)(576)$.
Then $B(\sigma)=(1,1,3,2,5,5,5)$.
By definition it is easy to see that $B$-code is a bijection from $\mathfrak{S}_n$ to $\mathcal{L}_n$.

\subsection{The set-valued statistic $\Lmicycl_1$}
Define the \emph{left-to-right minimum letters} statistic $\Lmil$, the set-valued \emph{left-to-right minimum places} statistic $\Lmip$ and the \emph{left to right minimun} statistic $\lrmin$ respectively by
$$\Lmil(\sigma):=\{\sigma_i:\,\sigma_i<\sigma_j\text{ for all }j<i\},$$
$$\Lmip(\sigma):=\{i:\,\sigma_i<\sigma_j\text{ for all }j<i\},$$
and
$$\lrmin(\sigma):= \# \Lmil(\sigma)\quad (\text{or }\# \Lmip(\sigma) ).$$
It is easy to see that $\Lmil(\sigma)=\Lmip(\sigma^{-1})$.

For an integer sequence $(\ell_1,\ldots, \ell_n)\in\mathcal{L}_n$, define $O((\ell_1,\ldots, \ell_n)):=\{i:\, \ell_i=1\}$, the set of indices with values $1$.
It turns out that $\Lmip$ and $\Lmil$ correspond to Lehmer code and $A$-code respectively.
The proof of the following lemma is directly by definition.

\begin{lem}
We have
\begin{equation}\label{eq:OA}
\Lmip(\sigma)=O(\Leh(\sigma)) \qquad \text{and} \qquad
\Lmil(\sigma)=O(A(\sigma)).
\end{equation}
\end{lem}

Hence it is natural to consider the statistic corresponding to $B$-code.
Define the set-valued statistic $\Lmicycl_1$ by
\begin{equation}~\label{eq:OB}
\Lmicycl_1(\sigma):=O(B(\sigma)).
\end{equation}
For example, $\Lmicycl_1(579328164)=O((1,1,3,3,1,6,2,6,3))=\{1,2,5\}$.

We have the following combinatorial interpretation of $\Lmicycl_1$, which explains the somewhat awkward notation, standing for the \emph{left-to-right minimum of the shifted cycle containing $1$}.

\begin{lem}
For $\sigma\in \mathfrak{S}_n$, write the the cycle containing $1$ in the way that $1$ is at the end of the cycle and denote the resulting cycle $\vec{c}$.
Then
$$\Lmicycl_1(\sigma) = \Lmil(\vec{c})$$
by regarding $\vec{c}$ as a word.
\end{lem}

\proof
By the definition of $B$-code, $b_i=1$ if and only if $i\in\vec{c}$ and all letters on the left of $i$ in $\vec{c}$ are
larger than $i$.
In other words, $i$ is a left-to-right minimum letter in $\vec{c}$.
Then we have $O(B(\sigma))=\Lmil(\vec{c})$.
\qed

\smallskip
For the running example, $\sigma=579328164=(1527)(394)(68)$ and the shifted cycle is $\vec{c}=(5271)$, hence
$\Lmicycl_1(579328164)=\Lmil(5271)=\{1,2,5\}$.

Also, $\Lmil$ and $\Lmicycl_1$ are related via $\phi$.

\begin{lem}\label{lem:afterphi1}
We have $$\Lmil(\sigma)=\Lmicycl_1(\phi(\sigma)).$$
\end{lem}
\proof
From \eqref{eq:OA}, \eqref{eq:OB} and the definition of $\phi$, we have
$$\Lmicycl_1(\phi(\sigma))=O(B(\phi(\sigma)))=O(A(\sigma))=\Lmil(\sigma).$$
\qed

%%%%%%%%%%%%%%%%%%%%%%%%%%%%%%%%%%%%%
\subsection{The set-valued statistics $\Inv$ and $\Sor$}
The goal of this subsection is to define and investigate the set-valued statistics $\Sor$ and $\Inv$.
First we need the concept of the \emph{induced set}.

\begin{defi}\label{defi:induce_set}
Given $(\ell_1, \ell_2,\ldots, \ell_n)\in \mathcal{L}_n$, define its
\emph{induced set} $\langle (\ell_1, \ell_2,\ldots, \ell_n) \rangle$
according to the following algorithm:
\begin{enumerate}
\item Set $S,U$ with the initial values $S=\{1,2,\ldots,n\}$ and
$U=\emptyset$.
\item For $i$ from $n$ down to $1$ do the followings:
    \begin{enumerate}
    \item let $\ell_i'$ be the $\ell_i$-th smallest element among $S$,
    \item add ordered pairs $(\ell_i',j)$ into $U$ for those $j\in S$ with $j>\ell_i'$. If there is no such $j$ then skip this step.
    \item delete $\ell_i'$ from $S$.
    \end{enumerate}
\item Define $\langle (\ell_1,\ell_2,\ldots, \ell_n)\rangle:=U.$
\end{enumerate}
\end{defi}

For example, for $(1,1,3,2,5,5,5)\in \mathcal{L}_7$ we have
\begin{center}
\begin{tabular}{c|c|c|l}
$i$ & $\ell_i'$ & $S$ & $U$ \\ \hline
  &   & 1,~2,~3,~4,~5,~6,~7 & $\emptyset$\\
7 & 5 & 1,~2,~3,~4,~*,~6,~7 & (5,6), (5,7) \\
6 & 6 & 1,~2,~3,~4,~*,~*,~7 & (5,6), (5,7), (6,7) \\
5 & 7 & 1,~2,~3,~4,~*,~*,~* & (5,6), (5,7), (6,7) \\
4 & 2 & 1,~*,~3,~4,~*,~*,~* & (5,6), (5,7), (6,7), (2,3), (2,4) \\
3 & 4 & 1,~*,~3,~*,~*,~*,~* & (5,6), (5,7), (6,7), (2,3), (2,4) \\
2 & 1 & *,~*,~3,~*,~*,~*,~* & (5,6), (5,7), (6,7), (2,3), (2,4), (1,3) \\
1 & 3 & *,~*,~*,~*,~*,~*,~* & (5,6), (5,7), (6,7), (2,3), (2,4), (1,3)
\end{tabular}
\end{center}
Thus $\langle (1,1,3,2,5,5,5)\rangle=\{(5,6), (5,7), (6,7), (2,3), (2,4), (1,3)\}$.

\medskip
Now we review the \emph{sorting index} $\sor$ of Petersen~\cite{Petersen_11}.
Given $\sigma$, decompose it uniquely into the product of transpositions $\sigma=(i_1j_1)(i_2j_2)\dots (i_kj_k)$ with $j_1<j_2<\dots <j_k$ and $i_r<j_r$ for $1\le r \le k$, and then define
$$\sor(\sigma):=\sum_{r=1}^{k}{(j_r-i_r)}.$$
For example, since $\sigma=2431765=(12)(24)(56)(57)$ we have $\sor(\sigma)=(2-1)+(4-2)+(6-5)+(7-5)=6$.
In other words, $\sor(\sigma)$ measures the total distance of the letters needed to move during the bubble-sorting process.
In this example, we have
$$2431756 \xrightarrow{(57)} 2431657 \xrightarrow{(67)} 2431567\xrightarrow{(24)} 2134567\xrightarrow{(12)} 1234567.$$

In~\cite{Chen_12} it is proved that $\inv(\sigma)=\sor(\phi(\sigma))$ and
\begin{equation}\label{eq:sor_Bcode}
\sor(\sigma)=\sum_{i=1}^n (i-b_i),
\end{equation}
which clarifies the relation between sorting index and the $B$-code $(b_1,b_2,\dots, b_n)$ of $\sigma$.
Observe that in the step (2)(b) of~\ref{defi:induce_set} we add exactly $(i-\ell_i)$ ordered pairs into $U$ for each $i$, hence from (\ref{eq:sor_Bcode}) it makes sense to define the set-valued statistic \emph{sorting set} $\Sor$ by
\begin{equation}\label{eq:Sor_Bcode}
\Sor(\sigma):=\langle  B(\sigma)\rangle.
\end{equation}

As for the $\Inv$, since $\inv(\sigma):=\#\{(i,j): i<j \text{ and } \sigma_i< \sigma_j\}$, it is natural to define the set-valued statistic \emph{inversion set} $\Inv$ by
$$\Inv(\sigma):=\{(i,j): i<j \text{ and } \sigma_i< \sigma_j\}.$$
Similar to the relation between $\Sor$ and $B$-code, we have the following:
\begin{pro}
For $\sigma\in\mathfrak{S}_n$, we have
\begin{equation}\label{eq:Inv_Acode}
\Inv(\sigma)=\langle A(\sigma)\rangle.
\end{equation}
\end{pro}
\proof
Let $\sigma=\sigma_1\cdots\sigma_n$ be the permutation with $A(\sigma)=(a_1,\ldots,a_n)$.
Since $A(\sigma)=\Leh(\sigma^{-1})$, hence $\sigma^{-1}_n=a_n$, $\sigma^{-1}_{n-1}=$ $(a_{n-1})$-th smallest element in $[n]\setminus\{\sigma^{-1}_n\}$.
In general, for $1< i <n$ we have
$$\sigma^{-1}_{n-i}= (a_{n-i})\text{-th smallest element in } [n]\setminus\{\sigma^{-1}_n, \ldots, \sigma^{-1}_{n-i+1}\}.$$

Thus $\sigma$ can be rebuilt from $A(\sigma)$ as follows.
At the initial stage there are $n$ vacancies from left to right.
For $i$ from $n$ down to $1$, we recursively put letter $i$ into the $a_i$-th vacancy from the left.
The resulting permutation is exactly $\sigma$.

For example, if $A(\sigma)=(1,1,3,2,5,5,5)$, then $\sigma$ can be recovered in the following way:
$$\_\,\_\,\_\,\_\,\_\,\_\,\_ \to \_\,\_\,\_\,\_7\_\,\_ \to \_\,\_\,\_\,\_76\_ \to \_\,\_\,\_\,\_765 \to \_\,\_4\_765 \to \_\,\_43765  \to 2\_43765 \to 2143765 $$

Now observe that in each step above, the position we choose for the letter $i$ is exactly $\ell_i'$ in Definition~\ref{defi:induce_set}.
In other words, we can synchronize the rebuilding of $\sigma$ and the construction of the induced set $\langle A(\sigma)\rangle$.
Moreover, note that the ordered pair $(\ell_i',j)$ is added to $U$ if and only if it is an inversion of $\sigma$, for the letter $i$ must be larger than $\sigma_j$.
Thus it must have $\Inv(\sigma)=\langle A(\sigma)\rangle$. \qed

\smallskip
And finally there is a set version of $\inv(\sigma)=\sor(\phi(\sigma))$:

\begin{lem}\label{lem:afterphi2}
We have $$\Inv(\sigma)=\Sor(\phi(\sigma)).$$
\end{lem}
\proof
From \eqref{eq:Sor_Bcode}, \eqref{eq:Inv_Acode} and the definition of $\phi$, we have
$$\Sor(\phi(\sigma))=\langle B(\phi(\sigma))\rangle = \langle A(\sigma) \rangle=\Inv(\sigma).$$
\qed

%%%%%%%%%%%%%%%%%%%%%%%%%%%%%%%%%%%%%%%%%%%%%%%%%%%%%%%%%%%%%%%%%%%%%%%%%%%%%%
%%%%%%%%%%%%%%%%%%%%%%%%%%%%%%%%%%%%%%%%%%%%%%%%%%%%%%%%%%%%%%%%%%%%%%%%%%%%%%
\section{Proof of Theorem 1.3}\label{sec:main_proof}
\noindent\emph{Proof of Theorem 1.3}.
(1) is obtained by combining Lemma~\ref{lem:afterphi1}, \ref{lem:afterphi2} and Theorem~\ref{thm:Foata_Han}.
The first statement of (2) is directly from (1) and in the following we look at the generating function.
For $n\geq 1$ let
$$F_n(q,x,y):=\sum_{\sigma\in\mathfrak{S}_n} q^{\sor(\sigma)}x^{\cyc(\sigma)}y^{\lrmincyc_1(\sigma)}.$$
It is clear that $F_1(q,x,y)=xy$. We claim that for $n\ge 2$ one has
$$F_n(q,x,y)=xy\prod_{r=2}^{n}(x+[r]_q+yq^{r-1}-1-q^{r-1}).$$
Let $t_{i\,j}$ denote the transposition $(ij)$ and let $\eta_1=1,\eta_2=1+t_{12}$ and
$$\eta_j=1+\sum_{i<j}t_{i\,j}$$
for $j\geq 3$.
Petersen~\cite{Petersen_11} showed that
\begin{equation}\label{eq:Sn_algebra}
\eta_1\eta_2\cdots\eta_n=\sum_{\sigma\in\mathfrak{S}_n}\sigma.
\end{equation}
Now define the linear map $\Theta:\mathfrak{S}_n\to\mathbb{Z}[q,x,y]$ by $\Theta(\sigma):=q^{\sor(\sigma)}x^{\cyc(\sigma)}y^{\lrmincyc_1(\sigma)}$.
Hence by \eqref{eq:Sn_algebra} it suffices to show that
$$\Theta(\eta_1\eta_2\cdots\eta_n)=xy\prod_{r=2}^{n}(x+[r]_q+yq^{r-1}-1-q^{r-1}).$$
It is easy to see that $\Theta(\eta_1)=xy$ and $\Theta(\eta_1\eta_2)=xy(x+yq)$.
Let $n\ge 3$. We proceed by induction. Suppose
$$\Theta(\eta_1\eta_2\cdots\eta_{n-1})=xy\prod_{r=2}^{n-1}(x+[r]_q+yq^{r-1}-1-q^{r-1}).$$

Take $\sigma = \sigma_1 \sigma_2\cdots\sigma_{n-1}
\in\mathfrak{S}_{n-1}$. It can be embedded in $\mathfrak{S}_n$ as
$\sigma = \sigma_1 \sigma_2\cdots\sigma_{n-1} n$. Let
$\sigma':=\sigma t_{i\,n}$ for some $1\leq i\leq n$ and it is clear
that $\sor(\sigma')=\sor(\sigma)+(n-i)$.

Assume that the cycle decomposition of $\sigma$ is $c_1\cdots c_m$ for some $m$ and $c_t$ contains the letter $i$.
Hence $\sigma'=c_1\cdots c_m(n)$ if $i=n$ and $c_1\cdots c'_t\cdots c_m$ if $i\neq n$, where $c'_t=(\ldots,i,n,\sigma(i),\ldots)$.
Thus
$$\cyc(\sigma')=
\begin{cases}
\cyc(\sigma)+1 & \text{if } i=n,\\
\cyc(\sigma) & \text{otherwise;}
\end{cases}
$$
and
$$\lrmincyc_1(\sigma')=
\begin{cases}
\lrmincyc_1(\sigma)+1 & \text{if } i=1,\\
\lrmincyc_1(\sigma) & \text{otherwise.}
\end{cases}
$$
So we have
\begin{eqnarray*}
\Theta(\sigma\cdot\eta_n) &=& \Theta(\sigma t_{n\,n}+\sigma t_{n-1\,n}+\cdots+\sigma t_{1\,n})\\
&=& \Theta(\sigma)(x+q+q^2+\cdots+q^{n-2}+yq^{n-1}),
\end{eqnarray*}
and therefore
\begin{eqnarray*}
\Theta(\eta_1\eta_2\cdots\eta_n)&=& \Theta\left(\sum_{\sigma\in\mathfrak{S}_n,\,\sigma(n)=n}\sigma\cdot\eta_n\right)\\
&=& \sum_{\sigma\in\mathfrak{S}_n,\,\sigma(n)=n}\Theta(\sigma\cdot\eta_n)\\
&=& (x+q+q^2+\cdots+q^{n-2}+yq^{n-1})\sum_{\sigma\in\mathfrak{S}_{n-1}}\Theta(\sigma)\\
&=& (x+[n]_q+yq^{n-1}-1-q^{n-1})\Theta(\eta_1\eta_2\cdots\eta_{n-1}).
\end{eqnarray*}
The proof is then completed by induction.
 \qed

%%%%%%%%%%%%%%%%%%%%%%%%%%%%%%%%%%%%%%%%%%%%%%%%%%%%%%%%%%%%%%%%%%%%%%%%%%%%%%
%%%%%%%%%%%%%%%%%%%%%%%%%%%%%%%%%%%%%%%%%%%%%%%%%%%%%%%%%%%%%%%%%%%%%%%%%%%%%%
\section{Set-valued joint equidistribution}\label{sec:triples}
In this section we seek for more set-valued statistics having the same joint distribution as $(\Rmil, \Lmil, \Lmap)$ and $(\Cyc, \Lmicycl_1, \Lmap)$.
In the introduction we have defined $\Rmil$ (right-to-left minimum letters) and $\Lmap$ (left-to-right maximum places) while in Section 2.2 $\Lmil$ (left-to-right minimum letters) and $\Lmip$ (left-to-right minimum places).
Similarly we can defined set-valued statistics $\Rmip, \Lmal, \Rmap$ and $\Rmal$ by
$$\Rmip(\sigma):= \{ i: \sigma_i < \sigma_j \text{ for all } j > i \},$$
$$\Lmal(\sigma):= \{ \sigma_i: \sigma_i > \sigma_j \text{ for all } j < i\},$$
$$\Rmap(\sigma):= \{i: i > \sigma_j \text{ for all } j > i \},$$
$$\Rmal(\sigma):= \{i: \sigma_i > \sigma_j \text{ for all } j > i \}.$$

The idea is quite simple: we look at the relations between these statistics by performing operations of ``inverse", ``complement", or ``reverse" on permutations.

For $\sigma=\sigma_1\sigma_2\dots \sigma_n \in \mathfrak{S}_n$, let $\sigma^{-1}$ denote its \emph{inverse},
$$\sigma^r:=(\sigma_n,\sigma_{n-1}, \dots ,\sigma_1)$$
its \emph{reverse} and
$$\sigma^c:=(n+1-\sigma_1, n+1-\sigma_2,\dots ,n+1-\sigma_n)$$
its \emph{complement}.
For example, if $\sigma=364152$, then $\sigma^{-1}=461352$, $\sigma^r=251463$ and $\sigma^c=413625$.
It is clear that the mappings $i,r,c: \mathfrak{S}_n\to \mathfrak{S}_n$, defined by $i(\sigma):=\sigma^{-1}$, $r(\sigma):=\sigma^r$ and $c(\sigma):=\sigma^c$, are bijections.

Given two set-valued statistics $\Stat_1$, $\Stat_2$ and one bijection $\chi: \mathfrak{S}_n\to \mathfrak{S}_n$, we say
$$\Stat_1\xrightarrow{\chi}\Stat_2$$
if $\Stat_1(\sigma)=\Stat_2(\chi(\sigma))$ for all $\sigma\in\mathfrak{S}_n$.
Also we define $\Stat^*:=\{n+1-i:\, i\in \Stat\}$ for a set-valued statistic $\Stat$, if applicable.

It turns out that these eight set-valued statistics are related via the mappings $i,r$ and $c$. The proof of the following proposition is straightforward and is omitted.

\begin{figure}[h]
  \includegraphics[width=2.5in]{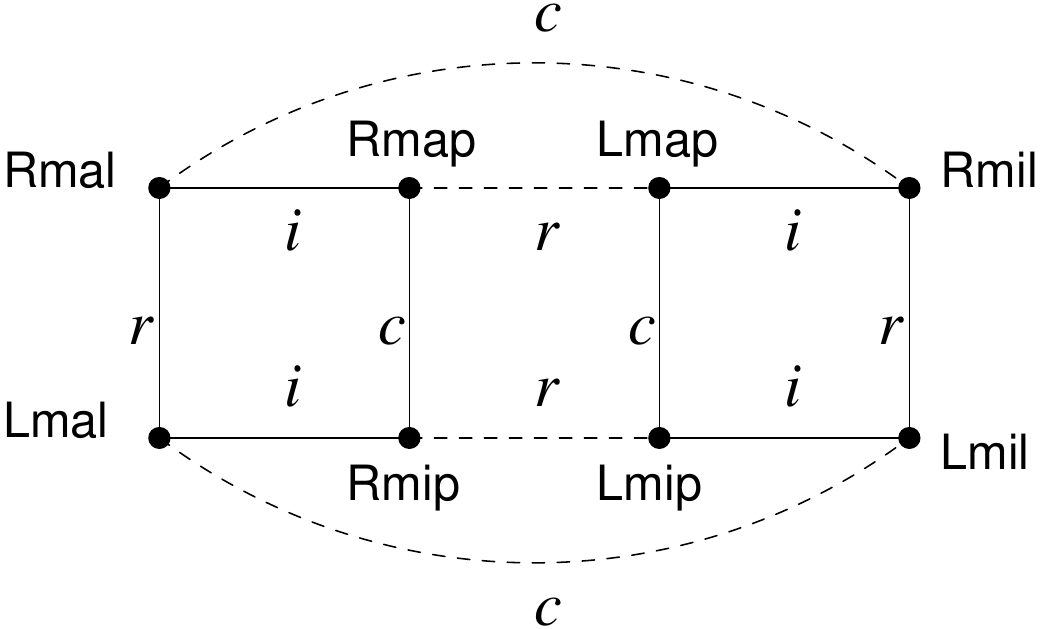}
  \caption{Relations between eight statistics \label{fig:8stat}}
\end{figure}

\begin{pro}~\label{prop:irc}
The relations between the eight set-valued statistics $\Lmil$, $\Lmip$,  $\Rmil$,  $\Rmip$,  $\Lmal$, $\Lmap$, $\Rmal$ and $\Rmap$ are illustrated by the graph (with edges solid or dotted and labeled by $i,r$ or $c$) in Figure~\ref{fig:8stat}.
\begin{enumerate}
\item A solid edge between $\Stat_1$ and $\Stat_2$ means $\Stat_1(\sigma)=\Stat_2(\chi(\sigma))$ and $\Stat_2(\sigma)=\Stat_1(\chi(\sigma))$ for $\sigma \in \mathfrak{S}_n$, with $\chi=i,r$ or $c$ as labeled.
\item A dotted edge between $\Stat_1$ and $\Stat_2$ means $\Stat_1(\sigma)=\Stat_2^*(\chi(\sigma))$ and $\Stat_2(\sigma)=\Stat_1^*(\chi(\sigma))$ for $\sigma \in \mathfrak{S}_n$, with $\chi=i,r$ or $c$ as labeled.
\end{enumerate}
\end{pro}

For examples, we have $\Rmap\xrightarrow{i}\Rmal$, $\Rmap\xrightarrow{r}\Lmap^*$, and $ \Rmap\xrightarrow{c}\Rmip$.

We come to the main result of this section.

\begin{thm} \label{thm:set_all}
The following triples of set-valued statistics have the same joint distribution over $\mathfrak{S}_n$.
A dash means the statistic is omitted.\\
\begin{tabular}{lll}
(1) $(\Rmil,   \Lmil,   \Lmap)$ & (2) $(\Cyc,    \Lmicycl_1, \Lmap)$ & (3) $(\Lmap,   \Lmip, \Rmil)$ \\
(4) $(\Rmip^*, \Rmap^*, \Lmal^*)$ & (5) $(\Lmal^*, \Rmal^*, \Rmip^*)$ & (6) $(\Lmil,\Rmil,\Rmap^*)$ \\
(7) $(\Rmal^*, \Lmal^*, \Lmip)$  & (8) $(\Rmap^*,\Rmip^*, \Lmil)$   & (9) $(\Lmip,   \Lmap, \Rmal^*)$ \\
(10) $(\Lmap, \text{-}  , \Cyc)$ & (11)$(\Lmicycl_1, \Cyc,  \text{-}  )$ & (12) $(\Cyc  ,  \text{-} , \Rmil)$  \\
(13) $(\Rmil, \text{-}  , \Cyc)$, &~  &~\\
\end{tabular}
\end{thm}

\proof By Proposition \ref{prop:irc}, we have
$$\begin{array}{cccc}
(1)\xrightarrow{i}(3), & (1)\xrightarrow{i\circ r \circ c}(4), & (1)\xrightarrow{r\circ c}(5), & (1)\xrightarrow{r}(6),\\
(1)\xrightarrow{c}(7), & (1)\xrightarrow{i\circ r}(8), & (1)\xrightarrow{i\circ c}(9).
\end{array}$$
Therefore, (3) to (9\textbf{}) are joint equidistributed with (1).
Moreover, by Theorem~\ref{thm:main}, we have
$$\begin{array}{ccccc}
(1)\xrightarrow{\phi}(2), & & (3)\xrightarrow{\phi}(10), & & (6)\xrightarrow{\phi}(11).
\end{array}$$
Finally, $(2)\xrightarrow{i}(12)$ and $(10)\xrightarrow{i}(13)$ follow from the fact that $\Cyc(\sigma)=\Cyc(\sigma^{-1})$, and the proof is completed. \qed

From the theorem we can read off the following many pairs of set-valued statistics which are symmetric and joint equidistributed.
Note that (1) includes those pairs of Foata and Han.

\begin{cor}
In each of the following items, the pairs of set-valued statistics are symmetric and joint equidistributed over $\mathfrak{S}_n$.
\begin{enumerate}
\item $(\Rmil,\Lmap)$, $(\Rmip^d,\Lmal^d)$, $(\Rmal^d,\Lmip)$,
$(\Rmap^d,\Lmil)$, $(\Cyc,\Rmil)$ and $(\Cyc,\Lmap)$
\item $(\Rmil,\Lmil)$, $(\Lmap,\Lmip)$, $(\Rmip^d,\Rmap^d)$,
$(\Rmal^d,\Lmal^d)$, and $(\Cyc,\Lmicycl_1)$.
\end{enumerate}
\end{cor}

%%%%%%%%%%%%%%%%%%%%%%%%%%%%%%%%%%%%%%%%%%%%%%%%%%%%%%%%%%%%%%%%%%%%%%%%%%%%%%
%%%%%%%%%%%%%%%%%%%%%%%%%%%%%%%%%%%%%%%%%%%%%%%%%%%%%%%%%%%%%%%%%%%%%%%%%%%%%%
\section{joint equidistributed quadruples}\label{sec:quadruples}
In this section we look at the ordinary number-valued statistics.
The goal is to find quadruples of statistics joint equidistributed with $(\inv, \rlmin, \lrmax, \lrmin)$ and $(\sor, \cyc, \lrmax, \lrmincyc_1)$.
Most of the materials in this section are well known.
Our contribution is to relate them with the $(\sor, \cyc, \lrmax, \lrmincyc_1)$ and derive the generating functions with respect to the first three statistics.

The statistics $\rlmin, \rlmax, \lrmin, \lrmax$ are defined in the obvious way.
We may take more familiar Mahonian statistics into consideration.
Let
$$\Des(\sigma):=\{i: \sigma_i>\sigma_{i+1}\}$$
be the \emph{descent set} of $\sigma$.
It is well known that the statistics \emph{major} $\maj$, \emph{inverse major} $\imaj$, \emph{reverse major} $\rmaj$, \emph{charge} $\chg$, and \emph{cocharge} $\cochg$, defined by
$$\maj(\sigma):=\sum_{i \in \Des} i,$$
$$\imaj(\sigma):=\maj(\sigma^{-1}), \qquad \rmaj(\sigma):=\maj(\sigma^{r\circ c}),$$
$$\chg(\sigma):=\sum_{i\in\Des(\sigma^{-1})}(n-i), \qquad \cochg(\sigma):=\sum_{i\notin\Des(\sigma^{-1})}(n-i),$$
are all Mahonian~\cite{Stanley_EC1}.

Again, given two statistics $\stat_1,\stat_2$ and one bijection $\chi: \mathfrak{S}_n \to \mathfrak{S}_n$, the notation
$$\stat_1\xrightarrow{\chi}\stat_2$$
means $\stat_1(\sigma)=\stat_2(\chi(\sigma))$ for all $\sigma\in\mathcal{S}_n$.
The equidistribution of $\maj$ and $\inv$ can be proved combinatorially from the celebrated \emph{fundamental bijection}
$\psi:\mathfrak{S}_n \to \mathfrak{S}_n$ of Foata~\cite{Foata_68}, that is, one has
$$\maj \xrightarrow{\psi} \inv.$$
Moreover, a close look of $\psi$ will also show that
$$ (\rlmax, \rlmin) \xrightarrow{\psi} (\rlmax, \rlmin).$$
We now perform the ``inverse-reverse-complement" trick on these statistics.

\begin{lem}
We have
$$\begin{array}{ccccccc}
(1)~\maj\xrightarrow{i}\imaj, & & (2)~\maj\xrightarrow{r\circ c}\rmaj, & &
(3)~\rmaj\xrightarrow{i}\chg, & & (4)~\chg\xrightarrow{r}\cochg.
\end{array}$$
\end{lem}
\proof
(1), (2) and (4) are obvious.
For (3), if $i\in\Des(\sigma^{r\circ c})$, then $(n+1)-\sigma^r_i < (n+1)-\sigma^r_{i+1}$ and thus $\sigma_{n-i}>\sigma_{n-i+1}$, which implies that $n-i \in\Des(\sigma)$.
Hence we have
$$ \rmaj(\sigma) = \maj(\sigma^{r\circ c}) = \sum_{i\in\Des(\sigma^{r\circ c})}i = \sum_{i\in\Des(\sigma)}(n-i) = \chg(\sigma^{-1}), $$
as desired. \qed

\smallskip
We are ready to state the main result of this section.

\begin{thm}~\label{thm:sta_all}
The following quadruple statistics are joint equidistributed over $\mathfrak{S}_n$:

\begin{tabular}{llll}
 (1) $(\inv,   \rlmin, \lrmin, \lrmax)$ & & (2) $(\sor,   \cyc,   \lrmincyc_1, \lrmax)$  \\
 (3) $(\inv,   \lrmax, \lrmin, \rlmin)$ & & (4) $(\inv,   \rlmin, \rlmax, \lrmax)$\\
 (5) $(\inv,   \lrmax, \rlmax, \rlmin)$ & & (6) $({n\choose 2}-\inv, \lrmin, \rlmin, \rlmax)$\\
 (7) $({n\choose 2}-\inv, \rlmax, \lrmax, \lrmin)$ & & (8) $({n\choose 2}-\inv, \rlmax, \rlmin, \lrmin)$\\
 (9) $({n\choose 2}-\inv, \lrmin, \lrmax, \rlmax)$ & & (10) $(\sor,   \lrmax, \lrmincyc_1, \cyc)$\\
 (11) $({n\choose 2}-\sor, \lrmincyc_1, \cyc,  \text{-} )$ & & (12) $({n\choose 2}-\sor,  \text{-},   \lrmax,  \lrmincyc_1)$ \\
 (13) $({n\choose 2}-\sor,  \text{-},   \cyc  ,  \lrmincyc_1)$ & & (14) $({n\choose 2}-\sor, \lrmincyc_1,   \lrmax  , \text{-})$ \\
 (15) $( \text{-}   , \cyc  ,  \text{-} , \rlmin)$ & & (16) $( \text{-}   , \rlmin,  \text{-} , \cyc)$\\
 (17) $(\maj,\rlmin,\rlmax, \text{-})$ & & (18) $(\imaj,\lrmax,\rlmax, \text{-})$\\
 (19) $(\rmaj,\lrmax,\lrmin, \text{-})$ & & (20) $(\chg,\rlmin,\lrmin, \text{-})$ \\
 (21) $(\cochg,\lrmin,\rlmin, \text{-})$, &~ &~\\
\end{tabular}

and the generating function with respect to the first three statistics in each quadruple is
$$F_n(q,x,y)=xy\prod_{r=2}^{n} (x+[r]_q+yq^{r-1}-1-q^{r-1}).$$
\end{thm}
\proof
Observe that $\inv(\sigma^{-1})=\inv(\sigma)$ and $\inv(\sigma^{r})=\inv(\sigma^{c})={n\choose 2}-\inv(\sigma)$ for
$\sigma \in \mathfrak{S}_n$.
The proof is done via the following mappings.
$$\begin{array}{ccccc}
(1)\xrightarrow{\phi}(2),  & (1)\xrightarrow{i}(3), & (1)\xrightarrow{i\circ r\circ c}(4), & (1)\xrightarrow{r\circ c}(5), & (1)\xrightarrow{r}(6), \\
(1)\xrightarrow{c}(7), & (1)\xrightarrow{i\circ r}(8), & (1)\xrightarrow{i\circ c}(9), & (3)\xrightarrow{\phi}(10), & (6)\xrightarrow{\phi}(11), \\
(6)\xrightarrow{\phi}(12), & (8)\xrightarrow{\phi}(13), & (9)\xrightarrow{\phi}(14), & (2)\xrightarrow{i}(15), & (10)\xrightarrow{i}(16)\\
(4)\xrightarrow{\psi}(17), &(17)\xrightarrow{i}(18), & (17)\xrightarrow{r\circ c}(19), & (19)\xrightarrow{i}(20), & (20)\xrightarrow{r}(21).
\end{array}$$
\qed

From the list we can read off many pairs of symmetric and joint equidistributed statistics.

\begin{cor}
In each of the following items, the pairs of statistics are symmetric and joint equidistributed over $\mathfrak{S}_n$.
\begin{enumerate}
\item $(\rlmax,\rlmin)$, $(\rlmax,\lrmax)$, $(\rlmin,\lrmin)$, $(\lrmax,\lrmin)$, $(\cyc,\lrmincyc_1)$, and $(\lrmincyc_1,\lrmax)$.
    The generating function with respect to each pair is $F(1,x,y)$.
\item $(\rlmax,\lrmin)$, $(\rlmin,\lrmax)$, $(\cyc,\lrmax)$, and $(\cyc,\rlmin)$.
\end{enumerate}
\end{cor}

%%%%%%%%%%%%%%%%%%%%%%%%%%%%%%%%%%%%%%%%%%%%%%%%%%%%%%%%%%%%%%%%%%%%%%%%%%%%%%
%%%%%%%%%%%%%%%%%%%%%%%%%%%%%%%%%%%%%%%%%%%%%%%%%%%%%%%%%%%%%%%%%%%%%%%%%%%%%%

\section{Concluding remarks}
In this short paper we generalize simultaneously Petersen and
Foata-Han's results to more than two statistics and find many
triples or quadruples of statistics having the same joint
distribution over $\mathfrak{S}_n$. In all quadruples, the first
statistic is Mahonian while the others are Stirling, and we then
read off many symmetric equidistributed pairs of Stirling
statistics. However, it is well known that there are pairs of
Mahonian statistics with a symmetric joint distribution as
well~\cite{Stanley_EC1}, for examples, $(\inv, \maj)$ and $(\maj,
\imaj)$ are two of them. Hence it would be interesting to generalize
our results further to include more Mahonian statistics.

On the other hand, the generating function obtained in Theorem~\ref{thm:main} only involves three of the four statistics, hence a natural question is to find a four-variable generating function including $\lrmax$ as well.
We leave them to the interested reader.

%%%%%%%%%%%%%%%%%%%%%%%%%%%%%%%%%%%%%%%%%%%%%%%%%%%%%%%%%%%%%%%%%%%%%%%%%%%%%%
%%%%%%%%%%%%%%%%%%%%%%%%%%%%%%%%%%%%%%%%%%%%%%%%%%%%%%%%%%%%%%%%%%%%%%%%%%%%%%
%%%%%%%%%%%%%%%%%%%%%%%%%%%%%%%%%%%%%%%%%%%%%%%%%%%%%%%%%%%%%%%%%%%%%%%%%%%%%%
%%%%%%%%%%%%%%%%%%%%%%%%%%%%%%%%%%%%%%%%%%%%%%%%%%%%%%%%%%%%%%%%%%%%%%%%%%
%%%%%%%%%%%%%%%%%%%%%%%%%%%%%%%%%%%%%%%%%%%%%%%%%%%%%%%%%%%%%%%%%%%%%%%%%%
%%%%%%%%%%%%%%%%%%%%%%%%%%%%%%%%%%%%%%%%%%%%%%%%%%%%%%%%%%%%%%%%%%%%%%%%%%
%%%%%%%%%%%%%%%%%%%%%%%%%%%%%%%%%%%%%%%%%%%%%%%%%%%%%%%%%%%%%%%%%%%%%%%%%%
\rm

%%%%%%%%%%%%%%%%%%%%%%%%%%%%%%%%%%%%%%%%%%%%%%%%%%%%%%%%%%%%%%%%%%%%%%%%%%
%%%%%%%%%%%%%%%%%%%%%%%%%%%%%%%%%%%%%%%%%%%%%%%%%%%%%%%%%%%%%%%%%%%%%%%%%%

\end{document}